\documentclass[12pt,a4paper]{article}
\usepackage{amsmath,amsfonts,amssymb,amscd}

\usepackage[a4paper]{geometry}

\usepackage{ulem}
\usepackage{xcolor}
\usepackage{graphics}

\usepackage{hyperref}

\def\ad{\operatorname{ad}}

\newcounter{th}
\def\t{\refstepcounter{th}{\bf \noindent{Theorem} \arabic{th}. }}

\newcounter{prop}
\def\prop{\refstepcounter{prop}{\bf \noindent{Proposition} \arabic{prop}. }}

\newcounter{lem}

\newcounter{de}
\def\de{\refstepcounter{de}{\bf \noindent{Definition} \arabic{de}. }}

\newcounter{ex}

\begin{document}

\begin{center}
    {\LARGE{\bf Double extension for commutative $n$-ary superalgebras with a skew-symmetric invariant form}}
\end{center}

\begin{center}
     E.G. Vishnyakova
\end{center}

\bigskip

\begin{abstract}
	
The method of double extension, introduced by A.~Medina and Ph.~Revoy, is a procedure which  decomposes a Lie algebra with an invariant symmetric form into elementary pieces. Such decompositions were developed for other algebras, for instance for Lie superalgebras and associative algebras, Filippov $n$-algebras and Jordan algebras. 

The aim of this note is to find a unified approach to such  decompositions using the derived bracket formalism. More precisely, we show that any commutative $n$-ary superalgebra with a skew-symmetric invariant form can be obtained inductively by taking orthogonal sums and generalized double extensions.

\end{abstract}

\bigskip

\section{Preliminaries}

Let $\mathfrak a = \mathfrak a_{\bar {0}} \oplus \mathfrak a_{\bar {1}}$ be a finite dimensional $\mathbb{Z}_2$-graded vector space over $\mathbb{K}$, where we assume for simplicity that $\mathbb{K} = \mathbb{R}$ or $\mathbb{C}$.  We denote by $\bar{a}\in \mathbb{Z}_2$ the parity of a homogeneous element $a\in \mathfrak a_{\bar {i}}$.  A bilinear form $(\,,)$ on $\mathfrak a$ is called
{\it even} if the corresponding  
linear map $\mathfrak a\otimes \mathfrak a \to \mathbb{K}$
is even.  A bilinear form is called {\it skew-symmetric} if
$$
(a,b)= -(-1)^{\bar{a}\bar{b}}(b,a)
$$
for any homogeneous elements $a,b\in \mathfrak a$.  
From now on we assume that $(\,,)$ is an even non-degenerate skew-symmetric  form on $\mathfrak a$. This is 
\begin{itemize}
	\item $(\,,)|_{\mathfrak a_{\bar {0}}\times \mathfrak a_{\bar {0}}}$ is a  non-degenerate skew-symmetric form, 
	\item $(\,,)|_{\mathfrak a_{\bar {1}}\times \mathfrak a_{\bar {1}}}$ is a  non-degenerate symmetric form, 
	\item $(\,,)|_{\mathfrak a_{\bar {0}}\times \mathfrak a_{\bar {1}}}=0$.
\end{itemize}

\medskip

\de \label{def invarian mult, symmetric}
$\bullet$ An {\it $n$-ary superalgebra} is a vector space $\mathfrak a$ together with an $n$-linear map 
$
\mathfrak a\times \cdots \times  \mathfrak a \longrightarrow \mathfrak a$. We denote this map by $(a_1,\ldots,a_n) \mapsto \{a_1,\ldots,a_n\}.$

$\bullet$ An $n$-ary superalgebra is called {\it commutative} if
\begin{equation}\label{eq symmetry}
	\{a_1,\ldots, a_i, a_{i+1}, \ldots, a_n \} = (-1)^{\bar{a_i}\bar{a}_{i+1}} \{a_1,\ldots, a_{i+1}, a_i, \ldots, a_n\}
\end{equation}
for all homogeneous $a_i, a_{i+1}\in \mathfrak a$.

\smallskip

$\bullet$ A commutative $n$-ary  superalgebra is called {\it invariant with respect to the form} $(\,,)$ if the following holds:
\begin{equation} \label{eq quadratic}
	(a_0, \{a_{1}, \ldots, a_{n}\}) = (-1)^{\bar{a}_0 \bar{a}_1}(a_1, \{a_0,a_2, \ldots, a_{n}\})
\end{equation}
for all homogeneous $a_{i}\in \mathfrak a$.

\medskip

We will write {\it a commutative invariant $n$-ary superalgebra} as a shorthand for {\it a commutative $n$-ary superalgebra that is invariant with respect to the form $(\,,)$}.
Let $\mathfrak a$ be a commutative $n$-ary superagebra. 

\medskip
\de\label{de n-ary subalgebra} $\bullet$ An {\it $n$-ary subalgebra} in $\mathfrak a$ is a vector subspace $\mathfrak b\subset\mathfrak a$  such that $\{\mathfrak b,\ldots,\mathfrak b\}\subset \mathfrak b$. 
 An {\it ideal} in $\mathfrak a$ is a vector subspace $\mathfrak i\subset \mathfrak a$  such that $\{\mathfrak a,\ldots,\mathfrak a,\mathfrak i\}\subset \mathfrak i$. 
\smallskip

$\bullet$ Let $\mathfrak a$ and $\mathfrak b$ be two $n$-ary algebras. An even linear map $\phi:\mathfrak a\to \mathfrak b$ is called a {\it homomorphism} of $n$-ary algebras if
$$
\phi(\{a_{1}, \ldots, a_{n}\}_{\mathfrak a}) = \{\phi(a_{1}), \ldots, \phi(a_{n})\}_{\mathfrak b}.
$$
The vector space $\ker\phi\subset \mathfrak a$ is an ideal in $\mathfrak a$. 

\smallskip

$\bullet$ A commutative $n$-ary superalgebra is called {\it simple} if it is not trivial one dimensional and it does not have any proper ideals.  

\smallskip

$\bullet$ An invariant commutative $n$-ary superalgebra is called {\it irreducible} if it is not a direct sum of two non-degenerate ideals. (An ideal $\mathfrak i$ is called {\it non-degenerate} if $(\,,)|_{\mathfrak i}$ is non-degenerate.)

\medskip

The main tool that we use in this paper is the derived braket construction. Let  $\mathfrak a$ be a $\mathbb Z_2$-graded vector space and $(\,,)$ be as above. We
denote by $S^n \mathfrak a$ the $n$-th symmetric power of $\mathfrak a$ and we put $S^*
\mathfrak a = \bigoplus\limits_n S^n \mathfrak a$. The superspace $S^* \mathfrak a$ possesses a natural
structure $[\,,]$ of a Poisson superalgebra that is defined in the following way:
$$
[x,y]: = (x,y), \quad x,y\in \mathfrak a;
$$
$$
[v,w_1\cdot w_2]: = [v,w_1]\cdot w_2 + (-1)^{vw_1} w_1\cdot [v,w_2],
$$
$$
[v,w] = -(-1)^{vw}[w,v],
$$
where $v,w, w_i$ are homogeneous elements in $S^* \mathfrak a$.  The super-Jacobi identity has the following form:
$$
[v,[w_1, w_2]] = [[v, w_1], w_2] + (-1)^{\bar{v} \bar{w_1}}[w_1,[v,w_2] ].
$$

Let us take any element $\mu \in S^{n+1} \mathfrak a$. Then we can define an n-ary superalgebra on $\mathfrak a$ in the following way:
\begin{equation}\label{eq der product}
	\{a_{1}, \ldots, a_{n}\}: = [a_1,[\ldots,[a_n,\mu]\ldots]], \,\,\, a_i\in \mathfrak a.
\end{equation}
Clearly, $\{a_{1}, \ldots, a_{n}\}\in \mathfrak a$. 
We will denote the corresponding superalgebra by $(\mathfrak a,\mu)$ and we will call the element $\mu$ the {\it derived
	potential} of $(\mathfrak a,\mu)$.  
The $n$-ary superalgebras $(\mathfrak a,\mu)$ are commutative and invariant with respect to the form $(\,,)$, see \cite{Ted Der Bracket homotopy} and \cite{Liza commutative algebras}.
We have the following observation (see \cite{Liza commutative algebras} for details).

\medskip

\prop\label{Main obseravation graded} {\it  Any commutative invariant $n$-ary superalgebras on $\mathfrak a$  can be obtained by construction
	(\ref{eq der product}). }

\medskip

\section{Double extension for invariant $n$-ary superalgebras}

 Let $\mathfrak{g}$ be an invariant commutative $n$-ary superalgebra and $\mu\in S^{n+1} \mathfrak{g}$ be its derived potential. Let $\mathfrak{h}$ be any commutative $n$-ary superalgebra with the multiplication $\nu \in S^{n} \mathfrak{h}^* \otimes \mathfrak{h}$. We can identify the vector spaces $S^{n} \mathfrak{h}^* \otimes \mathfrak{h}$  with the subspace $
  S^{n} \mathfrak{h}^* \cdot  \mathfrak{h} \subset S^*( \mathfrak{h} \oplus\mathfrak{h}^*)
  $
   and consider $\nu$ as a derived potential for the invariant superalgebra $\mathfrak{h} \oplus \mathfrak{h}^*$. The even non-degenerate skew-symmetric invariant form $(\,,)$ on $\mathfrak{h} \oplus \mathfrak{h}^*$ is given by:
   \begin{align*}
   (\alpha,x):= \alpha(x), \quad (x,\alpha):= -(-1)^{\bar{\alpha}\bar{x}} (\alpha,x),  
   \end{align*} 
   $\alpha\in \mathfrak{h}^*, \,\,x\in \mathfrak{h}$ are homogeneous elements.
   
   Let $\mathfrak{d} = \mathfrak{g} \oplus \mathfrak{h} \oplus \mathfrak{h}^*$ be a vector space with the non-degenerate skew-symmetric bilinear form $(\,,)$ that is the sum of the non-degenerate skew-symmetric bilinear forms on $\mathfrak{g}$ and $\mathfrak{h} \oplus \mathfrak{h}^*$. 

\medskip

\de\label{generalized DE for n-ary} The commutative invariant $n$-ary superalgebra $\mathfrak{d} = \mathfrak{g} \oplus \mathfrak{h} \oplus \mathfrak{h}^*$ with the derived potential 
\begin{equation}\label{eq der poten of double extension}
\begin{split}
L= \mu + \nu +\sum\limits_{i=1}^{n+1} \psi_i,\,\,\text{where} \quad \mu\in  S^{n+1} \mathfrak{g}, \\
 \nu \in S^{n} \mathfrak{h}^* \cdot\mathfrak{h}, \,\, \psi_i\in S^{i}\mathfrak{h}^* \cdot S^{n-i+1} \mathfrak{g}
\end{split}
\end{equation}
  is called a {\bf generalized double extension} of $\mathfrak{g}$ by $\mathfrak{h}$ via $\psi_i$, $i=1,\dots, n+1$. 

\medskip

The main observation of this section is:

\medskip

\t \label{theor decomposition general} {\it Assume that $\mathfrak{a}$ is an irreducible but not simple commutative invariant $n$-ary superalgebra. 
Then $\mathfrak{a}$ is isomorphic to a certain generalized double extension.}

\medskip

\noindent {\it Proof.} 
  Let us take a maximal non-trivial ideal $\mathfrak i$ of $\mathfrak{a}$. Clearly, $\mathfrak i^{\bot}$ is a minimal ideal in $\mathfrak{a}$. Furthermore, since $\mathfrak i^{\bot}\cap \mathfrak i$ is also an ideal and $\mathfrak{a}$ is irreducible, we see that $\mathfrak i^{\bot}\subset \mathfrak i$.  Therefore,  $\mathfrak i^{\bot}$ is a minimal isotropic ideal in $\mathfrak{a}$. Let us take a subspace $\mathfrak{h}$ in $\mathfrak{a}$ such that $\mathfrak{h}$ is isotropic, $\mathfrak{h} \cap  \mathfrak i^{\bot}=\{0\}$ and $(\,,)|_{\mathfrak{h}\oplus \mathfrak i^{\bot}}$ is non-degenerate. Since, $(\,,)|_{\mathfrak{h}\oplus \mathfrak i^{\bot}}$ is non-degenerate, we have $\mathfrak{a} = \mathfrak i\oplus \mathfrak{h}$.
 Consider the vector space $\mathfrak{w}= \mathfrak i^{\bot} \oplus \mathfrak{h}$.  We have the decompositions $\mathfrak{a} = \mathfrak{w} \oplus \mathfrak{w}^{\bot}$ and $\mathfrak i = \mathfrak i^{\bot} \oplus \mathfrak{w}^{\bot}$.

 Denote by $\lambda\in S^{n+1} \mathfrak{a}$ the derived potential of $\mathfrak{a}$.  Since, $\mathfrak{a} = \mathfrak{h} \oplus \mathfrak i^{\bot} \oplus \mathfrak{w}^{\bot}$, we have:
$$
S^{n+1} \mathfrak{a} = \bigoplus\limits_{i+j+k=n+1} S^{i} (\mathfrak{h})  \cdot S^{j} ( \mathfrak i^{\bot})  \cdot S^{k} (\mathfrak{w}^{\bot}).
$$ 
Therefore, we have the following decomposition of the derived potential: 
$$
\lambda = \sum\limits_{i+j+k=n+1} \lambda_{ijk},\,\,\,\text{where}\,\,\, \lambda_{ijk}\in S^{i} (\mathfrak{h})  \cdot S^{j} ( \mathfrak i^{\bot})  \cdot S^{k} (\mathfrak{w}^{\bot}).
$$ 
 Furthermore, for any $b\in \mathfrak i^{\bot}$ we have 
$$
[b,\lambda]\in \bigoplus\limits_{i+j+k=n+1} S^{i-1} (\mathfrak{h})  \cdot S^{j} (\mathfrak i^{\bot})  \cdot S^{k} (\mathfrak{w}^{\bot}).
$$
Since $\mathfrak i^{\bot}$ is an ideal we have 
$$
[a_1,\ldots, [a_{n-1},[b,\lambda]]]\in \mathfrak i^{\bot} \,\,\, \text{for} \,\,\, a_p\in \mathfrak{a}.
$$
Therefore, $\lambda_{ijk}$ can be non-trivial only in the following two cases:
\begin{itemize}
	\item $i=0$;
	\item $i=1$ and $k=0$.
\end{itemize}
In other words, we get
\begin{equation}\label{eq_form of der potential}
\lambda \in \Big(\mathfrak{h}  \cdot S^{n} (\mathfrak i^{\bot})\Big) \oplus \Big(\bigoplus\limits_{j+k=n+1} S^{j} ( \mathfrak i^{\bot})  \cdot S^{k} (\mathfrak{w}^{\bot}) \Big).
\end{equation}
We put 
$$
\begin{array}{l}
\mu:= \lambda_{0,0, n+1}\in S^{n+1} (\mathfrak{w}^{\bot});\quad
\nu:= \lambda_{1,n, 0}\in \mathfrak{h}  \cdot S^{n} ( \mathfrak i^{\bot});\\
\rule{0pt}{6mm}\psi_i:= \lambda_{0,i, n-i+1}\in S^{i} ( \mathfrak i^{\bot})  \cdot S^{n+1-i} (\mathfrak{w}^{\bot}), \,\,\, i=1,\ldots, n+1;\\
\rule{0pt}{6mm}\mathfrak{g}:= \mathfrak{w}^{\bot}. 
\end{array}
$$
We also can identify $\mathfrak i^{\bot} $ with $ \mathfrak{h}^*$.  We see that $(\mathfrak{a},\lambda)$ is a double extension of $(\mathfrak{g},\mu)$ by $(\mathfrak{h},\nu)$ via $\psi_i$. The result follows.$\Box$

\medskip

\prop\label{prop g is a factor} The $n$-ary superalgebra $(\mathfrak{g},\mu)$ is isomorphic to $\mathfrak i/\mathfrak i^{\bot}$. Indeed, assume that $a_i\in \mathfrak i$. Then
$$
[a_1,\ldots,[a_n,\nu]] = 0,\quad [a_1,\ldots,[a_n,\psi_i]]\in \mathfrak i^{\bot}. 
$$
Hence,
$$
\{\tilde a_1,\ldots, \tilde a_n\}_{\mathfrak i/\mathfrak i^{\bot}} = \{a_1,\ldots, a_n\}_{(\mathfrak{g},\mu)},
$$
where $\tilde a_i$ is the image of $a_i$ in $\mathfrak i/\mathfrak i^{\bot}$.

\medskip

\noindent{\bf Corollary 1.} {\it Assume that $\mathfrak{a}$ is irreducible and not simple, $\mathfrak i$ is a maximal non-trivial ideal and $\mathfrak{h}$ is an isotropic subalgebra in $\mathfrak{a}$ such that $\mathfrak{a} = \mathfrak i\oplus \mathfrak{h}$. Then $\mathfrak{a}$ is isomorphic to  a certain generalized double extension with  $\psi_{n}= \psi_{n+1}= 0$. In this case the generalized double extension is called {\bf double extension}.}

\medskip

\noindent {\bf Proof.} Consider (\ref{eq_form of der potential}) holds. 
 Let us take $x_i\in \mathfrak{h}$. The result follows from the following observations:
  $$
[x_1,\ldots,[x_n, \mu+\sum\limits_{i=1}^{n-1}\psi_{i}]] =0
$$
$$
 [x_1,\ldots,[x_n, \nu]]\in  \mathfrak{h}, \quad [x_1,\ldots,[x_n, \psi_{n}]]\in  \mathfrak{g},\quad [x_1,\ldots,[x_n, \psi_{n+1}]]\in  \mathfrak{h}^*.
$$
Since $\mathfrak{h}$ is a subalgebra, we have $\psi_{n}= \psi_{n+1} =0$. 
The proof is complete.$\Box$

\medskip

\noindent{\bf Corollary 2.} {\it 
	 Assume that $\mathfrak{a}$ is an irreducible but not simple skew-symmetric invariant $n$-ary algebra and  $\mathfrak i$ is a maximal non-trivial ideal of codimension $1$. Then $\mathfrak{a}$ is isomorphic to  a certain double extension with $\nu= \psi_i=0$ for all $i\ne 1$. 
	
	}

\medskip

\noindent {\bf Proof.} In this case the statement \ref{eq_form of der potential} has the following form:
$$
\lambda \in \mathfrak i^{\bot} \cdot S^{n} (\mathfrak{w}^{\bot}) \oplus S^{n+1} (\mathfrak{w}^{\bot}).\Box
$$

\subsection{\bf Lie algebras}

In this section we show that our definition coincides with the definitions given in \cite{MR} in case of Lie algebras. 

\medskip

\de A {\it derivation} of an $n$-ary algebra $(V,\mu)$ is a linear map $D:V\to V$ such that
$$
D(\{v_1,\ldots, v_n\}) = \sum_j \{v_1,\ldots, D(v_j),\ldots, v_n\}.
$$
We denote by $\mathrm{IDer} (\mu)$ the vector space of all derivations of the algebra $(V,\mu)$ preserving the form $(\,,)$. The following proposition was proven in \cite{Liza commutative algebras}:

\medskip

\prop\label{prop derivations} {\it Let us take any $w\in S^2(V)$ and $\mu \in S^{n+1}(V)$.
 We have:
	$$
	\mathrm{IDer} (\mu) =  \{w\in S^2(V)\,\, | \, \ad w(\mu) = 0  \}.
	$$
}

In \cite{MR} the following theorem was proven:

\medskip

\t\label{teor Medina, Revoy}{\bf [Medina, Revoy]} {\it  Let $\mathfrak{g}$ be a Lie algebra with an invariant symmetric form $B$. Let us take any  Lie algebra $\mathfrak{h}$ with a  homomorphism $\theta$ of $\mathfrak{h}$ to $\mathcal{D}er(\mathfrak{g})$ such that $\theta(\mathfrak{h})$ preserves the form $B$. Then $\mathfrak{d} = \mathfrak{h}^* \oplus \mathfrak{g} \oplus \mathfrak{h}$ is a Lie algebra with respect to the following multiplication:
$$
\begin{array}{l}
[(f_1, w_1,s_1), (f_2, w_2,s_2)] = \\ (\ad^*(s_2) f_1 - \ad^*(s_1) f_2 + \omega (w_1,w_2), [w_1,w_2] + \theta(s_1)(w_2) - \theta(s_2)(w_1), [s_1,s_2]  ).
\end{array}
$$ 
Here $\omega (w_1,w_2)(s) := B(\theta(s)w_1, w_2)$. The Lie algebra $\mathfrak{d}$ possesses an invariant form given by the sum of $B$ and the natural symmetric form on $\mathfrak{h}^* \oplus \mathfrak{h}$. (In \cite{MR} the algebra $\mathfrak{d}$ is called a {\bf double extension} of $\mathfrak{g}$ by $\mathfrak{h}$.)

Conversely, any Lie algebra with an invariant symmetric form can be obtained inductively by direct sums and double extensions.}

\medskip

We can simplify the definition of the double extension from \cite{MR} using the derived bracket construction. Let $\mathfrak{g}$ be a Lie algebra with an invariant symmetric form and $\mu$ be its derived potential. Consider a Lie algebra $\mathfrak{h}$ with the multiplication $\nu\in \bigwedge^2 \mathfrak{h}^* \otimes \mathfrak{h}$. We have the non-degenerate symmetric form on $\mathfrak{h}^* \oplus \mathfrak{h} \oplus \mathfrak{g}$: it is the sum of the form on $\mathfrak{g}$ and the natural non-degenerate symmetric form on $\mathfrak{h}^* \oplus \mathfrak{h}$.  As above we denote by $[\,,]$ the corresponding Poisson bracket on $S^*(\mathfrak{h}^* \oplus \mathfrak{h} \oplus \mathfrak{g})$.  (Again we can consider $\mathfrak{h}^* \oplus \mathfrak{h} \oplus \mathfrak{g}$ as a pure odd vector space with a non-degenerate skew-symmetric form.)

There is a one-to-one correspondence between elements $\psi\in \mathfrak{h}^* \otimes \bigwedge^2\mathfrak{g}$ and linear maps $\theta: \mathfrak{h} \to \bigwedge^2\mathfrak{g} \simeq \mathfrak{so}(\mathfrak{g})$. This correspondence is given by $\psi \longmapsto \theta_{\psi}$, where $\theta_{\psi}(x) = [x,\psi]$.

\medskip

\t\label{theor definitions of Double ext are equivalent}
Let $\mathfrak{d} = \mathfrak{h}^* \oplus \mathfrak{h} \oplus \mathfrak{g}$, and $\mathfrak{g}$, $\mathfrak{h}$ and $\theta = \theta_{\psi}$ be a double extension of  $\mathfrak{g}$ by $\mathfrak{h}$ via $\theta_{\psi}$ in the sense of Theorem \ref{teor Medina, Revoy}. Then in terms of Proposition \ref{Main obseravation graded} the Lie algebra $\mathfrak{d}$ has the derived potential 
$$
\mu +\nu+ \psi,
$$
and we have
\begin{equation}\label{eq mu + nu + phi = Medina}
	[\mu +\nu+ \psi,\mu +\nu+ \psi]=0.
\end{equation}

Conversely, if the condition (\ref{eq mu + nu + phi = Medina}) holds then 
the Lie algebra $\mathfrak{h}^* \oplus \mathfrak{h} \oplus \mathfrak{g}$ is a double extension of $\mathfrak{g}$ by $\mathfrak{h}$ via $\theta_{\psi}$ in the sense of Theorem \ref{teor Medina, Revoy}. 

\medskip

\noindent{\it Proof.} Assume that $\mathfrak{d}$ is a double extension of  $\mathfrak{g}$ by $\mathfrak{h}$ via $\theta_{\psi}$. An easy computation shows that the derived potential of $\mathfrak{d}$ is equal to $\mu +\nu+ \psi$. Since $\mathfrak{d}$ is a Lie algebra, we have (\ref{eq mu + nu + phi = Medina}).

Conversely, let us take Lie algebras $\mathfrak{g}$ and $\mathfrak{h}$ and an element $\psi$ as above such that (\ref{eq mu + nu + phi = Medina}) holds.  
 We need to show that $\theta = \theta_{\psi}$ is a homomorphism of $\mathfrak{h}$ to $\mathcal{D}er(\mathfrak{g})$ that preserves the form on $ \mathfrak{g}$.
 Since $\mathfrak{g}$ and $\mathfrak{h}\oplus \mathfrak{h}^*$ are Lie algebras we have $[\mu, \mu ] = 0$ and $[\nu,\nu] = 0$.  Therefore, we also have
 \begin{equation}\label{eq [L,L] =0}
 [\mu+\nu+ \psi, \mu+\nu+ \psi] = 2[\mu, \psi] + 2[\psi, \nu] + [\psi,\psi] = 0.
 \end{equation}
 Note that $[\mu, \nu] = 0$. Since $[\mu, \psi] \in \mathfrak{h}^* \otimes \bigwedge^3 \mathfrak{g}$ and $[\psi, \nu], [\psi,\psi] \in \bigwedge^2\mathfrak{h}^* \otimes \bigwedge^2 \mathfrak{g}$ we see that (\ref{eq [L,L] =0}) is equivalent to  
 \begin{equation}\label{eq 2[L,L] =0}
  [\mu, \psi] = 0, \quad 2[\psi, \nu] + [\psi,\psi] = 0
  \end{equation}

Let us show that 

\begin{itemize}
\item from $[\mu, \psi] = 0$ it follows  that $\theta(x)$ is a derivation of $\mathfrak{g}$ preser\-ving $(\,,)$;
\item from $2[\psi, \nu] + [\psi,\psi] = 0$ it follows that the map $\theta$ is a homomorphism from $\mathfrak{h}$ to $Der (\mathfrak{g})$
\end{itemize}

We have for any $x\in \mathfrak{h}$:
 $$
0=[x, [\mu, \psi]] = [[x,\mu], \psi] -[\mu, [x,\psi]] = -[\mu, [x,\psi]] = -[\mu, \theta(x)].
 $$
Now we apply Proposition \ref{prop derivations}. 

Let us study the second equation in (\ref{eq 2[L,L] =0}).  We need to show that the following holds:
$$
\theta([x,[y, \nu]]) = [\theta(x), \theta(y)] \,\, \text{for all}\,\,\, x,y\in \mathfrak{h}.
$$
Indeed,
$$
\theta([x,[y, \nu]]) = [\psi, [x,[y, \nu]]]  =  [[\psi,x],[y, \nu]] - [x,[[\psi,y], \nu]] + [x,[y, [\psi,\nu]]],\,\,  x,y \in \mathfrak{h}. 
$$
Notice that $[\psi,x] \in \bigwedge^2 \mathfrak{g}$. Therefore, $[[\psi,x],[y, \nu]] =0$. Similarly, $[x,[[\psi,y], \nu]] = 0$. Therefore,
$$
\theta([x,[y, \nu]])  =   [x,[y, [\psi,\nu]]],\,\,  x,y \in \mathfrak{h}. 
$$
Our statement follows from the following observation:
$$
 -\frac12 [x,[y, [\psi,\psi]]] = [x,[\psi, [y,\psi]]] = [[x,\psi], [y,\psi]] + 0 = [\theta(x), \theta(y)]. \quad \Box
$$

\medskip

Similar idea may be used to define a double extension for other types of quadratic algebras.

\subsection{\bf Solvable skew-symmetric invariant $n$-ary algebras}

Define by induction the following subalgebras in $\mathfrak a$.
\begin{align*}
	\mathfrak a^{(1)}:= \{\mathfrak a,\ldots,\mathfrak a\},\quad \mathfrak a^{(2)}:= \{\mathfrak a^{(1)},\ldots,\mathfrak a^{(1)}\},\cdots 
\end{align*}
By induction, we see that  $\mathfrak a^{(k+1)}\subset \mathfrak a^{(k)}$ and more precisely that $\mathfrak a^{(k+1)}$ is an ideal in $\mathfrak a^{(k)}$. We call $n$-ary superalgebra $\mathfrak a$ {\it solvable}  if there exists an integer $K$ such that $\mathfrak a^{(K)}=\{0\}$.

\medskip
\prop\label{prop subalgebra and hom are solvable} {\it Subalgebras and homomorphic images of a solvable commutative $n$-ary superalgebra are solvable. }

\medskip

\noindent {\it Proof.} Let $\mathfrak b$ be a subalgebra in a solvable $n$-ary superalgebra $\mathfrak a$. Then we see 
$$
\mathfrak b^{(1)} = \{\mathfrak b,\ldots,\mathfrak b\}  \subset \{\mathfrak a,\ldots,\mathfrak a\}= \mathfrak a^{(1)}.
$$
By induction we have $\mathfrak b^{(k)} \subset \mathfrak a^{(k)}$ for all $k$. Hence, $\mathfrak b$ is solvable.

Let $\phi:\mathfrak a\to \mathfrak h$ be a homomorphism. Denote by $\mathfrak b$ the image of $\phi$. Then, 
$$
\mathfrak b^{(1)} = \{\phi(\mathfrak a),\ldots,\phi(\mathfrak a)\} = \phi (\{\mathfrak a,\ldots,\mathfrak a\}) = \phi(\mathfrak a^{(1)}).
$$ 
Again by induction we have $\mathfrak b^{(k)} = \phi(\mathfrak a^{(k)})$.  Hence, $\mathfrak b$ is solvable. The proof is complete.$\Box$

\medskip

	In particular, all solvable skew-symmetric $n$-ary algebras can be obtained inductively by such double extensions.

Clearly, $\mathfrak a^{(1)} = \{\mathfrak a,\ldots,\mathfrak a\}$ is an ideal in $\mathfrak a$ and any proper subspace $\mathfrak i$ in $\mathfrak a$ such that $\mathfrak a^{(1)}\subset \mathfrak i$ is an ideal in $\mathfrak a$. Hence, there exists a maximal idea of codimension $1$. The result follows from Theorem \ref{theor decomposition general}.

The second observation follows from $\mathfrak g \simeq \mathfrak i/ \mathfrak i^{\bot}$, hence $\mathfrak g$ is solvable.

\subsection*{Acknowledgments}

Research of E.~Vishnyakova was supported by AFR-grant, University of Luxembourg and MPI Bonn.

\noindent{\it Elizaveta Vishnyakova}

\noindent {Max Planck Institute for Mathematics Bonn}

\noindent{\emph{E-mail address:} 	\verb"VishnyakovaE@googlemail.com"}


\begin{thebibliography}{99}

\bibitem[ABB2]{ABB reductive even} {\it  Albuquerque, H., Barreiro, E., Benayadi, S.}  Quadratic Lie superalgebras with a reductive even part. J. Pure Appl. Algebra 213 (2009), no. 5, 724-731.


\bibitem[BBB]{generalized double ext}{\it Bajo, I., Benayadi, S., Bordemann, M.} Generalized double extension and descriptions of quadratic Lie superalgebras
arXiv: math-ph/0712.0228 (2007)


\bibitem[BB]{Extension for Lie super} {\it Benamor, H., Benayadi,  S.} Double extension of quadratic Lie superalgebras. 
Comm. Algebra 27 (1999), no. 1, 67-88. 


\bibitem[Bord]{Bordemann} {\it Bordemann, M.} Nondegenerate invariant bilinear forms on nonassociative algebras. 
Acta Math. Univ. Comenian. (N.S.) 66 (1997), no. 2, 151-201. 


\bibitem[FF]{Figueroa-O'Farrill}{\it Figueroa-O'Farrill, J.} Metric Lie n-algebras and double extensions.  	arXiv:0806.3534





\bibitem[HW]{HW}{\it Hanlon P., Wachs M. L.} On Lie k-algebras, Adv. in Math. 113 (1995) 206-236.


\bibitem[KoSch1]{Kosman1} {\it Kosmann-Schwarzbach, Y.} Derived brackets. Letters in Mathematical Physics 69: 61-87, 2004.


\bibitem[MR]{MR} {\it  Medina, A., Revoy, Ph.} Algèbres de Lie et produit scalaire invariant. (French) [Lie algebras and invariant scalar products] Ann. Sci. École Norm. Sup. (4) 18 (1985), no. 3, 553-561.


\bibitem[V]{Liza commutative algebras} {\it Vishnyakova E.} Commutative $n$-ary superalgebras with an invariant skew-symmetric form.  arXiv:1409.4342.

\bibitem[Vor1]{Ted Der Bracket homotopy} {\it Voronov, Th.} Higher derived brackets and homotopy algebras. J. Pure Appl. Algebra 202 (2005), no. 1-3, 133-153.



\end{thebibliography}
\end{document}